\documentclass[11pt,a4paper]{article}
\usepackage[utf8]{inputenc}
\usepackage[english]{babel}
\usepackage[T1]{fontenc}
 \pdfoutput=1
\usepackage{amscd,amssymb,amsmath,amsthm,amstext,amsfonts}
\usepackage{stmaryrd}
\usepackage[pdftex]{graphicx}
\usepackage{fancybox}
\usepackage[square,numbers]{natbib}
\usepackage{booktabs}
\usepackage{hyperref}
\usepackage{algorithm}
\usepackage{algpseudocode}
\usepackage{enumerate}  
\usepackage[capitalise,nameinlink]{cleveref}
\crefformat{equation}{\textup{#2(#1)#3}}
\crefname{figure}{Figure}{Figures}
\crefformat{subfigure}{Figure #2#1#3}
\usepackage{comment}
\usepackage{mathtools}
\usepackage{tikz,pgfplots}
\pgfplotsset{compat=newest}
\usepackage{pgfplotstable}
\usepackage{caption}
\usepackage{subcaption}
\pgfplotsset{compat=newest}
\usepackage{multirow}
\usepackage{authblk} 
\usepackage[normalem]{ulem}

\theoremstyle{definition}

\theoremstyle{remark}

\numberwithin{theorem}{section}
\numberwithin{equation}{section}
\numberwithin{table}{section}
\numberwithin{figure}{section}

\allowdisplaybreaks

\author[1]{Moataz M. Alghamdi}
\affil[1]{Applied Mathematics and Computational Sciences (AMCS), King Abdullah University of Science and Technology. Thuwal, 23955-6900, Kingdom of Saudi Arabia}

\author[2]{Fleurianne Bertrand}
\affil[2]{Faculty of Electrical Engineering, Mathematics and Computer Science, University of Twente. Zilverling P.O. Box 217, 7500 AE Enschede, The Netherlands}

\author[1,4]{Daniele Boffi}

\author[3]{Francesca Bonizzoni}
\affil[3]{Institute of Mathematics, University of Augsburg, Universit\"atsstr.~12a, 86159 Augsburg, Germany}

\author[1]{Abdul Halim}

\author[1]{Gopal Priyadarshi}

\affil[4]{Department of Mathematics ``F. Casorati'', University of Pavia. Via Ferrata 1, 27100 Pavia, Italy}

\title{On the matching of eigensolutions to parametric partial differential equations}

\date{}

\begin{document}

\maketitle

\begin{abstract}
In this paper a novel numerical approximation of parametric eigenvalue problems is presented. We motivate our study with the analysis of a POD reduced order model for a simple one dimensional example. In particular, we introduce a new algorithm capable to track the matching of eigenvalues when the parameters vary.
\end{abstract}

\vspace{1cm}

\section{Introduction}

The study of parametric eigenvalue problems arising from partial differential equations with multidimensional (possibly stochastic) parameter space, is still the object of very limited research.
Starting from the pioneer work presented in~\cite{AndreevSchwab2012} it is  apparent that the analysis of parametric eigenvalue problems cannot be simply considered as a generalization of the theory developed for parametric/stochastic \emph{source} partial differential equations~\cite{BonizzoniNobile2020,BonizzoniNobileKressner,BonizzoniNobile2014,BonizzoniBuffaNobile13,BonizzoniNobile12,NobileTemponeWebster,BaeckNobileTamelliniTempone}. Indeed, parametric eigenvalue problems lack a fundamental regularity condition needed for the analysis of parametric source problems; this is the consequence of possible eigenvalue crossings occurring when the parameters vary. When a crossing occurs, clearly the eigenvalues involved in the crossing are not smooth functions of the parameters and the corresponding eigenspaces are not even continuous if the eigenvalues are sorted by their magnitude. The question addressed in this paper concerns the matching of the eigenvalues across their intersections so that a new sorting of the eigenmodes can be introduced that restores the smoothness of eigenvalues and eigenspaces with respect to the parameters.

A reduced basis approximation of an isolated eigenmode has been presented and analyzed in~\cite{FumagalliManzoniParoliniVerani}, while in~\cite{HorgerWohlmuthDickopf} the reduced basis model approach is applied to the simultaneous approximation of multiple eigenvalues. The latter reference can be considered as the state of the art in reduced order modeling for eigenvalue problems.

A reduced order model for the approximation of eigenvalue problems was considered in~\cite{Abdul} and an algorithm for tracking the matching of the eigenvalues is under development~\cite{tbd}. The latter takes inspiration from reduced order model techniques for the parametric-in-frequency Helmholtz equation~\cite{BonizzoniPradoveraRuggeri,BonizzoniNobilePerugiaPradovera20a,BonizzoniNobilePerugiaPradovera20b,BonizzoniNobilePerugia18,BonizzoniPradovera,Nobile-Pradovera}.

The aim of this paper is twofold.
On one side, we provide the reader with a convincing example of the necessity of tracking the matching of the eigenvalues for different parameters values. Indeed, the lack of prior knowledge of the behavior of the eigenvalues, in terms of dependence on the parameters and of their possible crossings, may lead to unexpected results.
On the other side, we introduce a greedy algorithm that can be used to successfully match the eigenmodes and we describe some of its properties.

In Section~\ref{se:ps} we describe our abstract problem. In Section~\ref{se:motiv} we present a one dimensional example from which it is clear how crucial is to detect the crossings of eigenvalues, and finally in Section~\ref{se:alg} we introduce our matching algorithm.

\section{Problem setting}
\label{se:ps}

Let $(H,(\bullet,\bullet)_H)$ and $(V,(\bullet,\bullet)_V)$ be Hilbert spaces such that $V\subset H\simeq H'\subset V'$ gives a standard Hilbert triplet and $V$ is compact subset of $H$. 
Moreover, let $\mathcal M\subset \mathbb R^P$ be a $P$-dimensional parametric domain, with $P\geq 1$, and $a,b\colon V\times V\times \mathcal M\rightarrow \mathbb R$
two parameter-dependent bilinear forms such that, for all $\mu\in\mathcal M$, $a(\bullet,\bullet;\mu)$ is symmetric and coercive, namely, there exist a positive constant $\alpha$ such that
\begin{equation}
\label{eq:b_property}
\aligned
&a(v,v;\mu)\geq \alpha \|v\|^2_V &&\forall\,v\in V\\
&a(w,v;\mu)=a(v,w;\mu) &&\forall\, w,v\in V
\endaligned
\end{equation}
and $b(\bullet,\bullet;\mu)$ is equivalent to the scalar product of $H$, namely, there exist positive constants $c_b,\, C_b$ such that
\begin{equation}
\label{eq:b_form}
c_b(w,v)_H\leq b(w,v)\leq C_b(w,v)_H\quad\forall\, w,v\in V.
\end{equation}
Given a window of values $[\lambda_{min},\lambda_{max}]\subset \mathbb R_+$ we are interested in the following parametric eigenvalue model problem: for each $\mu\in\mathcal{M}$, find eigenvalues $\lambda(\mu)\in[\lambda_{min},\lambda_{max}]$ and non-vanishing eigenfunctions $u(\mu)\in V$ such that, for all $v\in V$ it holds 
\begin{equation}
\label{eq:model_problem}
a(u(\mu),v;\mu)=\lambda(\mu) b(u(\mu),v;\mu).
\end{equation}

\section{A motivating example}
\label{se:motiv}

Let $\mathcal{M}$ be the interval $[-0.9,0.9]$ and consider the following $\mu$-dependent boundary value problem, with $\mu\in\mathcal M$:
\begin{equation}
\label{eq:PDE_example}
\left\{
\begin{array}{ll}
-\operatorname{div}(A(\mu)\nabla u(\mu))=\lambda(\mu)u(\mu)&  \textrm{ in }\Omega=(0,1)^2\\
u(\mu)=0& \textrm{ on }\partial\Omega 
\end{array}
\right.
\end{equation}
where the diffusion $A(\mu)\in\mathbb R^{2\times 2}$ is given by  the diagonal matrix
\begin{equation*}
A(\mu)\coloneqq 
\begin{pmatrix}
1& 0\\
0 &1+\mu
\end{pmatrix}.
\end{equation*}
The weak formulation of~\eqref{eq:PDE_example} reads: for all $\mu\in\mathcal M$, find $(\lambda(\mu),u(\mu))\in \mathbb R_+\times H^1_0(\Omega)$, with $u(\mu)$ non vanishing, such that, for all $v\in H^1_0(\Omega)$ it holds
\begin{equation}
\label{eq:weak_example}
\int_\Omega (A(\mu)\nabla u(\mu))\cdot\nabla v\, dx
=\lambda(\mu)\int_\Omega u(\mu) v\, dx.
\end{equation}
Problem~\eqref{eq:weak_example} is a particular case of the general problem~\eqref{eq:model_problem} when choosing the spaces $V=H^1_0(\Omega)$, $H=L^2(\Omega)$, (equipped with the natural inner products $(\bullet,\bullet)_{H^1}$ and  $(\bullet,\bullet)_{L^2}$) and the bilinear forms
\begin{equation}
\label{eq:a_b_example}
\begin{aligned}
a(w,v;\mu)&\coloneqq\int_\Omega (A(\mu)\nabla w)\cdot\nabla v\, dx,\\
b(w,v;\mu)&\coloneqq\int_\Omega w v\, dx.
\end{aligned}
\end{equation}
Notice that $a(\bullet,\bullet;\mu)$ is symmetric and coercive, with coercivity constant $\alpha=(1+C_P^2)^{-1}$, $C_P$ being the Poincar\'e constant, and $b(\bullet,\bullet;\mu)$ coincides with the $L^2$-inner product, i.e., the chain of inequalities~\eqref{eq:b_form} is a chain of equalities with constants $c_b=C_b=1$.

\begin{figure}
	\begin{subfigure}[b]{0.45\textwidth}
		\centering
		\includegraphics[width=\textwidth]{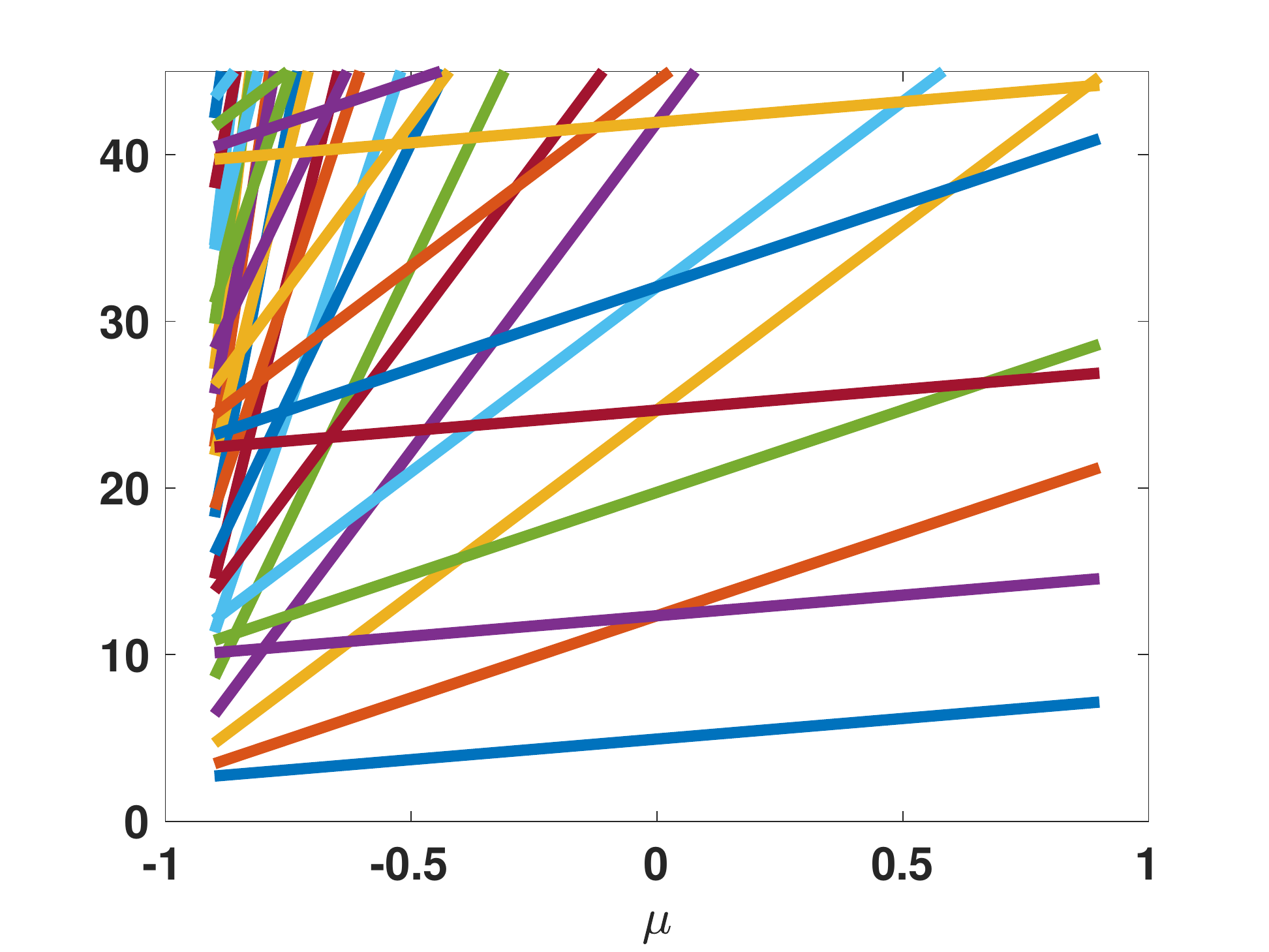}
		\caption{}
		\label{fig:exact_eig}
	\end{subfigure}
	\hfill
	\begin{subfigure}[b]{0.45\textwidth}
		\centering
		\includegraphics[width=\textwidth]{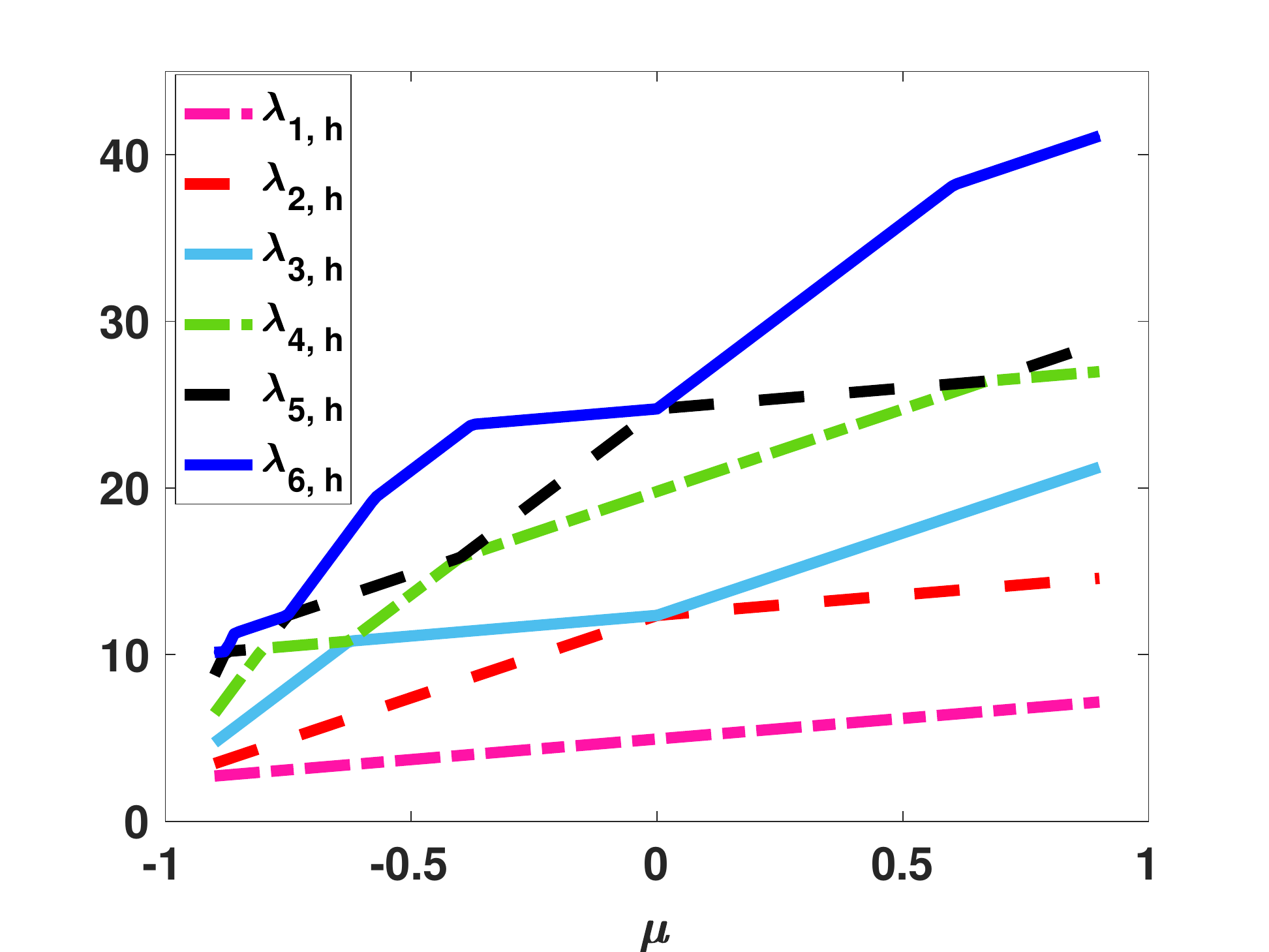}
		\caption{}
		\label{fig:six_eig}
	\end{subfigure}
	\caption{\ref{fig:exact_eig}: Exact eigenvalues $\lambda_{n,m}(\mu)$ given by formula~\eqref{eq:exact_eigensol}, for $\mu\in\mathcal M$.  \ref{fig:six_eig}: First six eigenvalues $\lambda_{1,h},\dots,\lambda_{6,h}$ computed by the FEM and sorted according to their magnitutde, for $\mu\in\mathcal M$.}
\end{figure}

The analytical eigensolutions to~\eqref{eq:PDE_example} can be explicitly computed by separation of variables, and they are given by:
\begin{equation}
\label{eq:exact_eigensol}
\begin{array}{ll}
\lambda_{n,m}(\mu)&=\frac{\pi^2}{4}(m^2+(1+\mu)n^2)\\
u_{n,m}&=\cos\left(\frac{m\pi}{2}x\right)\cos\left(\frac{n\pi}{2}y\right)
\end{array}
\quad\forall\, m,n\in\mathbb N.
\end{equation}
In particular, we underline that the eigenfunctions $\{u_{n,m}\}_{n,m\in\mathbb N}$ are independent of the parameter $\mu\in\mathcal{M}$, in contrast to the eigenvalues $\{\lambda(\mu)_{n,m}\}_{n,m\in\mathbb N}$. This property makes the problem particularly simple and suitable for our preliminary considerations. Figure~\ref{fig:exact_eig} shows the exact eigenvalues corresponding to the formula in Equation~\eqref{eq:exact_eigensol}, while Figure~\ref{fig:six_eig} shows the first six approximating eigenvalues sorted according to their magnitude and denoted $\lambda_{1,h},\dots,\lambda_{6,h}$.

\subsection{Reduced basis approximation of the first eigenvalue}

We are interested in computing an approximation to the first eigenpair $(\lambda_1(\mu),u_1(\mu))$ as $\mu$ varies in the parametric interval $\mathcal M$. Consider the uniform decomposition of $\mathcal M$
\begin{equation*}
\mathcal M_{T}=\{-0.9+(j-1)\Delta\mu,\, j=1,\ldots,T\} 
=\{\mu_j,\, j=1,\ldots,T\}
\end{equation*}
with $\Delta\mu=0.1$ and $T=19$. On a given regular (fine) mesh of $\Omega$, we compute the eigensolutions $\{(\lambda_{1,h}^{(j)},u^{(j)}_{1,h}),\, j=1,\ldots,T\}$ corresponding to $\mathcal{M}_T$ via the piecewise linear finite element method (FEM). We collect the $T$ computed eigenfunctions into the snapshot matrix $S_1=[u^{(1)}_{1,h}|\cdots|u^{(T)}_{1,h}]\in\mathbb R^{N_h\times T}$, where $N_h$ denotes the number of degrees of freedom. By performing the singular value decomposition (SVD), we derive the following representation of the snapshot matrix:
\begin{equation*}
S_1=U\Sigma Z^T,
\end{equation*}
where $U\in\mathbb R^{N_h\times N_h}$, $Z\in\mathbb R^{T\times T}$ are unitary matrices and $\Sigma\in\mathbb R^{N_h\times T}$ is a rectangular diagonal matrix.

\begin{figure}
	\centering
	\includegraphics[width=0.6\textwidth]{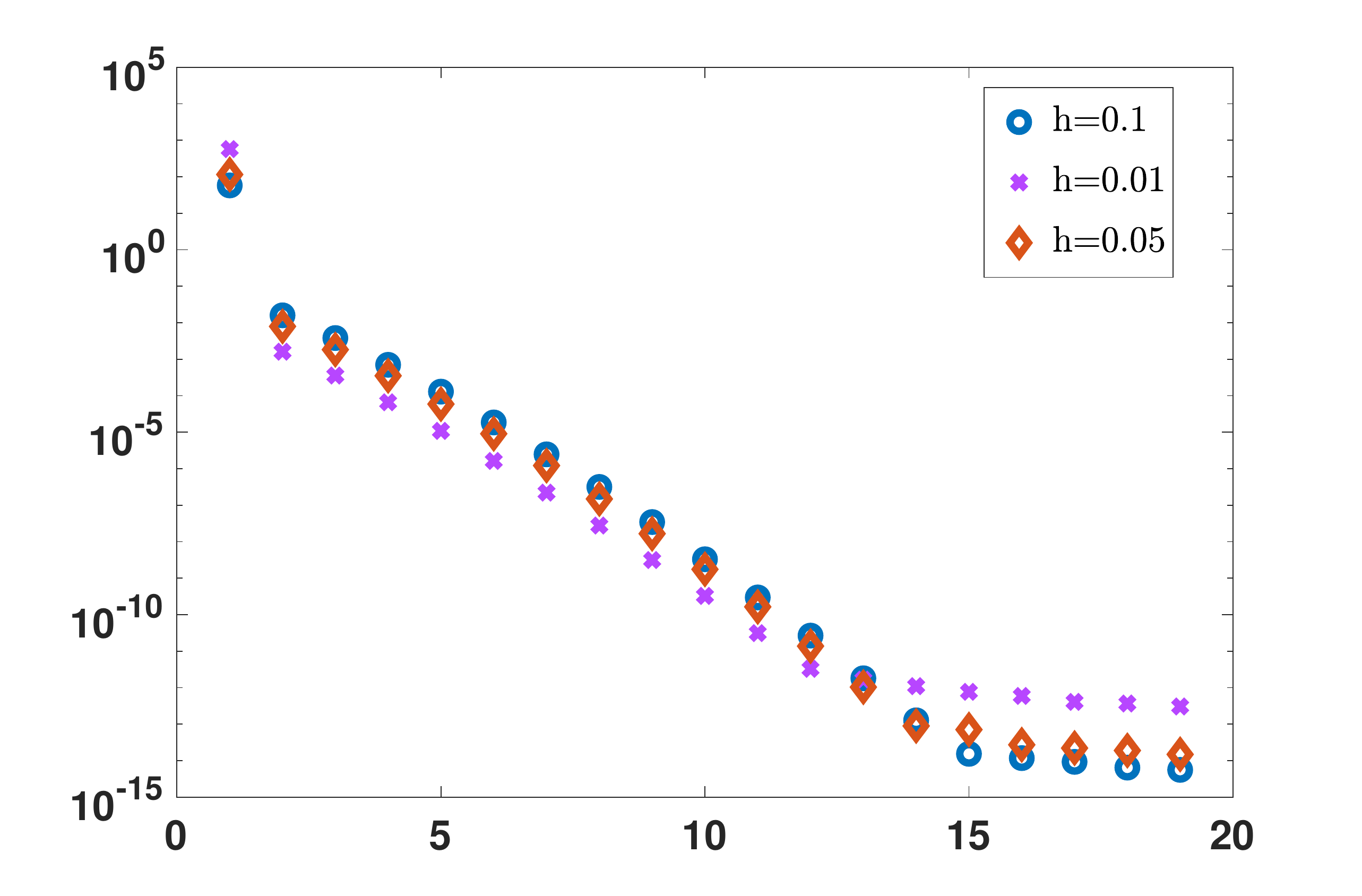}
	\caption{Singular values of the snapshot matrix $S_1$ for successively refined meshes of $\Omega$ with maximum diameter denoted as $h$. }
	\label{fig:singular_values}
\end{figure}

Theoretically we expect $S$ to have rank one because the first eigenvalue is well separated by the others and the first eigenspace is independent of $\mu$. Numerically, we observe that the first singular value is well separated by the others even if it is not the only non vanishing one. Indeed, for any fixed (relatively large) tolerance $tol<1.e-3$, several singular vectors $N_{tol}$ will be considered in the truncated SVD expansion of $S_1$ (see Figure~\ref{fig:singular_values}). 
For $tol>1.e-1$, only the first singular vector will be considered in the singular value decomposition, namely, $N_{tol}=1$. The reduced basis proper orthogonal decomposition (RB-POD) approximation to the first eigenpair $(\lambda_1(\mu),u_1(\mu))$ of~\eqref{eq:weak_example} is obtained by projection onto the one-dimensional space spanned by the first singular eigenvector. Looking at the results summarized in Table~\ref{tab:ROM_1}, we note that the first RB-POD eigenvalue is a good approximation of the first FE eigenvalue $\lambda_1(\mu)$, for $\mu\in\{-0.75,-0.25,0.25,0.75\}$. Slightly better approximations are also obtained for $N_{tol}=2$ (see Table~\ref{tab:ROM_2}).

\begin{table}[]
	\centering
	\begin{tabular}{ |c|c|c|c|} 
		\hline
		$h$ & $\mu$ & FEM based first eigenvalue & RB-POD based first eigenvalue \\
		\hline
		\multirow{3}{4em}{0.1} & -0.75 &  3.09172930& 3.09178369\\ 
		& -0.25 &4.32853369		 & 4.32853489	\\ 
		& 0.25 &  5.56526834&   5.56528610\\ 
		& 0.75 &  6.80197424&6.80203730 \\ 
		\hline
		\multirow{3}{4em}{0.05} & -0.75 & 3.08606437&3.08607518	\\ 
		& -0.25 & 4.32052203	&   4.32052233\\ 
		& 0.25 & 5.55496589	&  5.55496949\\ 
		& 0.75 & 6.78940395	& 6.78941665	\\ 
		\hline
		\multirow{3}{4em}{0.01} & -0.75 & 3.08432204 & 3.08432252\\ 
		& -0.25 & 4.31805168& 4.31805169	\\ 
		& 0.25 & 5.55178071 & 5.55178087\\ 
		& 0.75 & 6.78550948& 6.78551005\\ 
		\hline
	\end{tabular}
	\caption{Comparison between the FE and RB-POD approximation to the first eigenpair $(\lambda_1(\mu),u_1(\mu))$ of problem~\eqref{eq:weak_example}, using a one-dimensional RB space, spanned by the first singular vector of $S_1$.}
	\label{tab:ROM_1}
\end{table}

\begin{table}[]
	\centering
	\begin{tabular}{ |c|c|c|c|} 
		\hline
		$h$ & $\mu$ & FEM based first eigenvalue & RB-POD based first eigenvalue \\
		\hline
		\multirow{3}{4em}{0.1} & -0.75 &  3.09172930& 3.09172950                
		\\ 
		& -0.25 &4.32853369              & 4.32853469\\ 
		& 0.25 &  5.56526834&   5.56526837\\ 
		& 0.75 &  6.80197424&6.80197627 \\ 
		\hline
		\multirow{3}{4em}{0.05} & -0.75 & 3.08606437                            &3.08606442                              \\ 
		& -0.25 & 4.32052203    &   4.32052230\\ 
		& 0.25 & 5.55496589     &  5.55496590\\ 
		& 0.75 & 6.78940395     & 6.78940447\\ 
		\hline
		\multirow{3}{4em}{0.01} & -0.75 & 3.08432204 & 3.08432204               
		\\ 
		& -0.25 & 4.31805168& 4.31805169        \\ 
		& 0.25 & 5.55178071 & 5.55178071\\ 
		& 0.75 & 6.78550948& 6.78550950\\ 
		\hline
	\end{tabular}
	\caption{Comparison between the FE and RB-POD approximation to the first eigenpair $(\lambda_1(\mu),u_1(\mu))$ of problem~\eqref{eq:weak_example}, using a two-dimensional RB space, spanned by the first two singular vectors of $S_1$.}
	\label{tab:ROM_2}
\end{table}

\subsection{Reduced basis approximation of the third eigenvalue}

\begin{figure}
	\centering
	\includegraphics[width=0.6\textwidth]{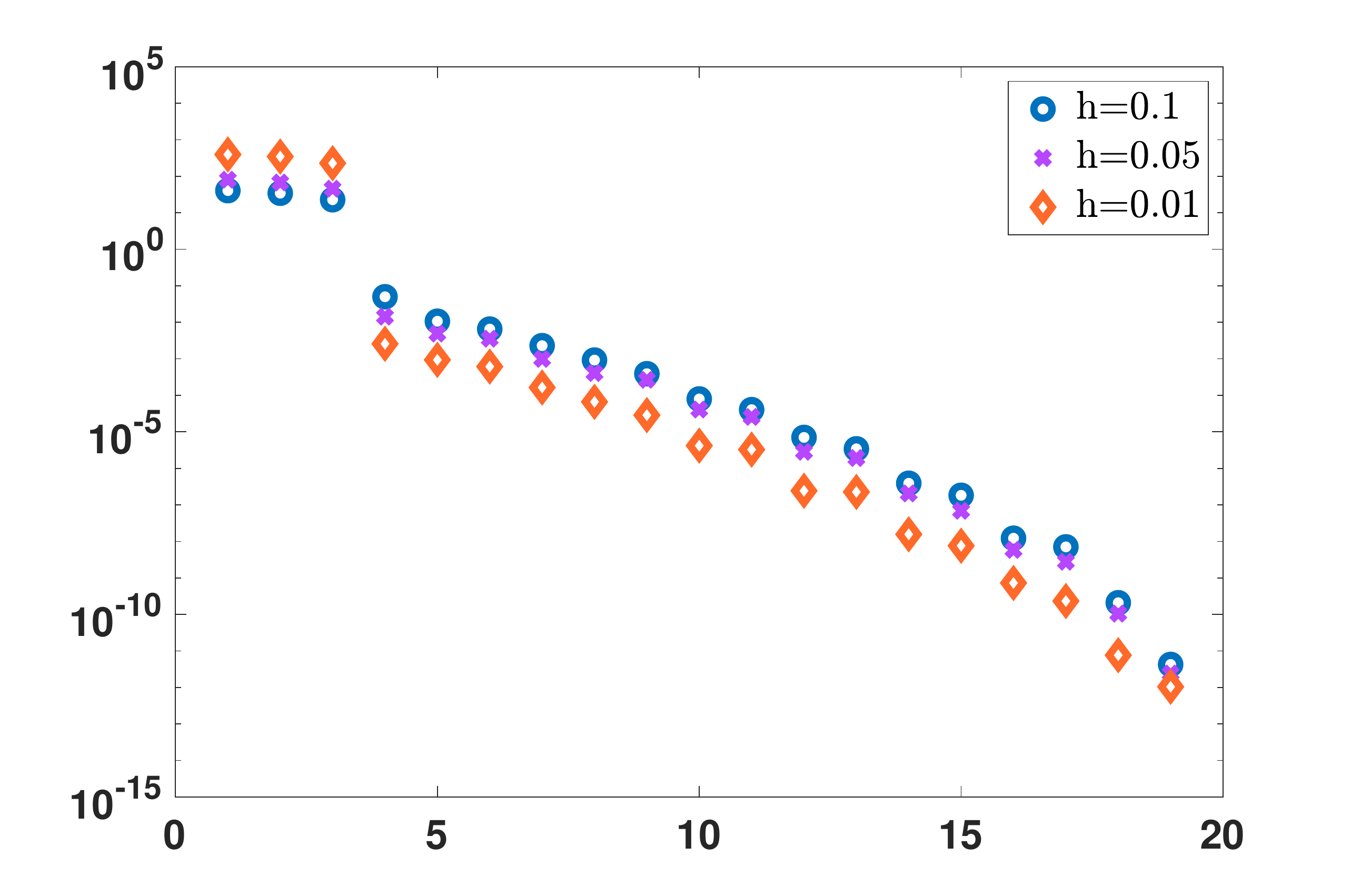}
	\caption{Singular values of the snapshot matrix $S_3$.}
	\label{fig:singular_values2}
\end{figure}

We now follow the same strategy as before, with the aim of approximating the third eigenpair $(\lambda_3(\mu),u_3(\mu))$ of problem~\eqref{eq:weak_example}, for $\mu\in\mathcal M$. Denote by $S_3$ the snapshot matrix $S_3=[u^{(1)}_{3,h}|\cdots|u^{(T)}_{3,h}]\in\mathbb R^{N_h\times T}$ collecting the third eigenfunctions $u^{(j)}_{3,h}$ with $\mu_j\in\mathcal M_T$. 
Theoretically, we expect $S_3$ to have rank 3, because of two eigenvalue crossings (see Figure~\ref{fig:six_eig}).
The singular values of $S_3$ are depicted in Figure~\ref{fig:singular_values2}, and the approximation results for $N_{tol}=3$ are summarized in Table~\ref{tab:ROM_3}. Even though the results might look satisfactory, it is important to observe that the numbers reported in the last column of Table~\ref{tab:ROM_3} correspond to the second eigenvalue of the $3\times 3$ reduced model.

Actually, in this case we know the exact solution and, after careful inspection, it was possible to realize that the approximation of the solution we are interested in, corresponds to the second eigenvalue of the $3\times 3$ system. This comes from the fact that the three element of the reduced basis correspond to the three eigenvalues associated with the three modes belonging to the third eigenmode. More precisely, looking at Figure~\ref{fig:six_eig}, the curve corresponding to $\lambda_{3,h}$, is made of three straight pieces and that's the reason why we are expecting the rank of the snapshot matrix to be equal to three. If we now isolate from the figure of the exact values~\ref{fig:exact_eig} the three straight lines corresponding to the three selected eigenfunctions, then we see that the curve we are interested in is always the one related to the second eigenfunction out of those three.

We can deduce that, in general, it is essential to know some information about the structure of the exact solution. In particular, a fundamental question that needs to be addressed is how to match computed eigenvalues for different values of the parameter $\mu$.

\begin{table}[]
	\centering
	\begin{tabular}{ |c|c|c|c|c|c|} 
		\hline
		$h$ & $\mu$  & FEM based third eigenvalue & RB-POD based third eigenvalue \\
		\hline
		\multirow{3}{4em}{0.1} & -0.75  &8.14338931                             
		
		&8.14352843                                                     \\ 
		& -0.25 & 11.78888922            &11.78893305\\ 
		& 0.25  &  14.89196477 &   14.89205929\\ 
		& 0.75  &  19.85303433           &19.85335723 \\ 
		\hline
		\multirow{3}{4em}{0.05} & -0.75 & 8.05008647                            
		&8.05010991             
		\\ 
		& -0.25 & 11.73700317   &   11.73701402\\ 
		& 0.25 & 14.82575077    & 14.82577826\\ 
		& 0.75 & 19.76686171                    & 19.76695884   \\ 
		\hline
		\multirow{3}{4em}{0.01} & -0.75 & 8.02024667                            
		& 8.02024755           
		\\ 
		& -0.25 & 11.72081569   & 11.72081613   \\ 
		& 0.25 & 14.80523805 & 14.80523912\\ 
		& 0.75 & 19.74028492& 19.74028861       \\ 
		\hline
	\end{tabular}
	\caption{Comparison between the FE and RB-POD approximation to the third eigenpair $(\lambda_3(\mu),u_3(\mu))$ of problem~\eqref{eq:weak_example}, when the RB space is spanned by the first three singular vectors of $S_3$.}
	\label{tab:ROM_3}
\end{table}

\section{The matching of eigenvalues}
\label{se:alg}

The motivating example of the previous section demonstrates that the reduced order techniques available in the literature are inappropriate for tracking the eigenpair solutions to the parameter-dependent eigenvalue problem~\eqref{eq:model_problem}. This section describes the crucial ingredient of a novel algorithm under development, which is able to overcome this obstacle.

The algorithm performs an a priori matching between two sets of eigensolutions. Given two values of the parameter $\mu_i,\,\mu_k\in\mathcal{M}$, and the corresponding set of eigenpairs
\begin{equation*}
\{(\lambda_j(\mu_i),u_j(\mu_i))\}_{j= 1}^{m_i}\qquad
\{(\lambda_\ell(\mu_k),u_\ell(\mu_k))\}_{\ell= 1}^{m_k}
\end{equation*}
for each $j=1,\ldots,m_i$, we want to find the value of $\ell\in\{1,\ldots,m_k\}$ such that $\lambda_j(\mu_i)$ and $\lambda_\ell(\mu_k)$ belong to the same eigenvalue curve $\lambda^\star\colon\mathcal{M}\rightarrow[\lambda_{min},\lambda_{max}]$, i.e.,
\begin{equation}
\label{eq:lambda_curve}
\lambda^\star(\mu_i)=\lambda_j(\mu_i),\qquad
\lambda^\star(\mu_k)=\lambda_\ell(\mu_k).
\end{equation}
Notice that this problem may have no solutions if some curve $\lambda^\star$ is entering or exiting the window $[\lambda_{min},\lambda_{max}]$ in the interval identified by $\mu_i$ and $\mu_k$. For the same reason, in general, $m_i$ may be different from $m_k$.

For this purpose, we adapt to our setting the following idea from~\cite{Nobile-Pradovera}. Under the assumption $m_i=m_k=m$, we construct the cost matrix $D\in\mathbb R^{m\times m}$
\begin{equation}
\label{eq:cost_matrix}
D_{j,\ell}^{i,k} = |\lambda_j(\mu_i)-\lambda_\ell(\mu_k)|
+ w \min(\|u_j(\mu_i)-u_\ell(\mu_k)\|,\|u_j(\mu_i)+u_\ell(\mu_k)\|)
\end{equation}
where $w$ is a suitable positive weight. We aim at finding one value per row and one value per column of $D$ so that the sum of the selected matrix entries is minimized. In other words, we look for a permutation $\boldsymbol{\sigma}=(\sigma_1,\ldots,\sigma_m)\colon\{1,\ldots,m\}\rightarrow\{1,\ldots,m\}$ such that $\lambda_j(\mu_i)$ and $\lambda_{\sigma_j}(\mu_k)$ belong to the same eigenvalue curve in the sense of equation~\eqref{eq:lambda_curve}, for $j=1,\ldots,m$. This is an optimization problem for which various solutions methods are available; for instance, a quite convenient solution strategy involves the use of the Hungarian algorithm.

We make a couple of observations.
\begin{itemize}
	\item 
	Each entry of the cost matrix~\eqref{eq:cost_matrix} has two ingredients: the first measures the distance between the two sets of eigenvalues, and the second measures the distance between the two sets of eigenfunctions. The weight $w$ express the relative importance of the second term with respect to the first one. Even though one might be tempted to consider the first term, only, i.e., taking $w=0$, in the majority of the cases this might lead to the wrong matching.
	
	\item
	In applications we often get $m_i\neq m_k$, leading to a rectangular cost matrix $D$. Typically, this happens when an eigenvalue curve $\lambda^\star(\mu)$, $\mu\in\mathcal M$, attains values that are out of the window of interest $[\lambda_{min},\lambda_{max}]$. In this situation, the cost matrix is rectangular, and the output of the Hungarian matrix is a permutation matching $m=\min\{m_i,m_k\}$ eigenpairs.  
\end{itemize}

Preliminary computations show that the a priori matching performs generally well with some exceptions.
In particular, in some cases the matching strategy described above might fail, delivering the wrong eigenpair matching and, moreover, it might not be able to deal with clusters of eigenvalues, namely, when two or more eigenvalues are close to each other, even if not multiple. For these reasons, a novel a posteriori matching strategy is under development which is able to resolve these issues.
Starting from an initial discretization $\mathcal{M}_0$ of the parameter set $\mathcal{M}$, and employing the a posteriori indicator, we will be able to build up a greedy algorithm that selects the areas of $\mathcal{M}$ where refinement is needed, delivering a problem-adapted discretization $\mathcal{M}_1$ of $\mathcal{M}$. 

\section*{Acknowledgements}
The work of F.\ Bertrand, D.\ Boffi, and A.\ Halim was supported by the Competitive Research Grants Program CRG2020 ``Synthetic data-driven model reduction methods for modal analysis'' awarded by the King Abdullah University of Science and Technology (KAUST).
D.\ Boffi is member of the INdAM Research group GNCS and his research is partially supported by IMATI/CNR and by PRIN/MIUR.
F.\ Bonizzoni is member of the INdAM Research group GNCS and her work is part of a project that has received funding from the European Research Council ERC under the European Union's Horizon 2020 research and innovation program (Grant agreement No.~865751).

\end{document}